\documentclass[letterpaper, 10 pt, conference]{ieeeconf}  % Comment this line out if you need a4paper

\IEEEoverridecommandlockouts                              % This command is only needed if 
\usepackage{algorithm}
\usepackage{algpseudocode}
\usepackage{graphicx}
\usepackage{tabularx}
\usepackage{algpseudocode}
\usepackage{amssymb}
\usepackage{cite}
\usepackage{amsmath}
\usepackage[utf8]{inputenc} % allow utf-8 input
\usepackage[T1]{fontenc}    % use 8-bit T1 fonts
\usepackage{hyperref}       % hyperlinks
\usepackage{url}            % simple URL typesetting
\usepackage{booktabs}       % professional-quality tables
\usepackage{amsfonts}       % blackboard math symbols
\usepackage{nicefrac}       % compact symbols for 1/2, etc.
\usepackage{microtype}      % microtypography
\usepackage{xcolor}         % colors
\usepackage{mathtools}
\usepackage{algorithm}
\usepackage{caption}
\usepackage{subcaption}
\usepackage{color}

\providecommand{\com[2]}{\begin{tt}[#1 #2]\end{tt}}

\title{On the Worst-Case Analysis of Cyclic Coordinate-Wise Algorithms
\\ on Smooth Convex Functions}

\author{Yassine Kamri, Julien M. Hendrickx, François Glineur
   \\
  ICTEAM, UCLouvain\\
  Louvain-La-Neuve, 1348, Belgium\\
  \texttt{\{yassine.kamri,julien.hendrickx,francois.glineur\}@uclouvain.be} \\
}

\begin{document}

\maketitle

\begin{abstract}
 We propose a unifying framework for the automated computer-assisted worst-case analysis of cyclic block coordinate algorithms in the unconstrained smooth convex optimization setup.
We compute exact worst-case bounds for the cyclic coordinate descent and the alternating minimization algorithms over the class of smooth convex functions, and provide sublinear upper and lower bounds on the worst-case rate for the standard class of functions with coordinate-wise Lipschitz gradients. We obtain in particular a new upper bound for cyclic coordinate descent that outperforms the best available ones by an order of magnitude. We also demonstrate the flexibility of our approach by providing new numerical bounds using simpler and more natural assumptions than those normally made for the analysis of block coordinate algorithms. Finally, we provide numerical evidence for the fact that a standard scheme that provably accelerates random coordinate descent to a $O(1/k^2)$ complexity is actually inefficient when used in a (deterministic) cyclic algorithm.

\end{abstract}

\section{Introduction}
% In this paper, we study the worst-case behaviour of cyclic block-coordinate descent algorithm (BCD) on the set of block-wise smooth functions. Coordinate descent like algorithm have been extensively studied because of their practical performance for large-scale optimization problems that arises for instance in machine learning. Worst-case bounds for cyclic BCD are derived in \cite{Beck}. Random coordinate descent and accelerated versions of coordinate descent are also extensively studied notably in \cite{Nesterov} for random BCD and in \cite{richtarik} for acceleration. Random BCD and acceleration are out of the scope of this article, we refer the reader to \cite{wright} for an detailed survey on coordinate descent like algorithms.
Large-scale optimization problems are the cornerstone of many engineering applications, such as in machine learning or signal processing. With the widespread availability of data, the scale of some of these problems is constantly increasing, to the point where standard full-gradient optimization methods are often becoming computationally too expensive. Fortunately, many of these problems possess a structure that allows the use of \emph{partial gradient} methods. An important subclass of these methods are the block coordinated descent algorithms, which only need access to a subset of gradient coordinates at each iteration. These can generally be separated in three main categories depending on how the blocks of coordinates are selected and updated \cite{Shi2016,wright}: (i) Gauss-Southwell methods greedily select the coordinates that lead to the largest improvements (i.e the coordinates with largest gradient norm), (ii) randomized methods select coordinates according to a probability distribution and (iii) cyclic methods update the coordinates in a cyclic predefined order. Although greedy methods can exhibit good performances, their update rule typically requires having access to full gradients. Hence randomized and cyclic methods have been more heavily used and studied.

The theoretical convergence analysis of random coordinate descents has proven easier than the analysis of their deterministic counterparts. Sampling coordinates with replacement from a suitable probability distribution implies indeed that the expectation of each coordinate step is the full gradient, making the analysis mostly similar to that of full gradient descent. Consequently, many random coordinate descent algorithms with theoretical guarantees have been proposed for convex optimization problems, including accelerated and proximal variants for a variety of probability distributions \cite{Nesterov,Lin,richtarik,Diakonikolas,Allen,Hanzely,Stich}. However, these performances are only guaranteed in expectation or with high probability, and the sampling technique can be computationally costly. Hence there is also an interest for cyclic block coordinate methods, which appear simpler and more efficient to implement in practice, but much harder to analyse. The difficulty appears to reside in establishing a link between the block of coordinates at each step and the full gradient. Some convergence results are already known for cyclic coordinate descent, but they are often obtained under quite restrictive assumptions such as the isotonicity of the gradient \cite{Saha}, or with worst-case bounds that are relatively conservative and under initial conditions that are not entirely standard for unconstrained convex optimization \cite{Beck}. Better convergence results exist for the specific case of quadratic optimization problems \cite{Hong,Li,Wang,Ching,Mert,goujaud}.

In this paper, we propose an alternative analysis based on Performance Estimation Problems (PEP). The underlying idea of PEP, which is to compute performance guarantees for first-order methods thanks to semidefinite programming (SDP) originated in \cite{Drori}. It was further developed in \cite{Taylor} where the authors used convex interpolation to derive tight results and provide examples of worst-case functions for various types of first-order methods. A related approach was proposed in \cite{Lessard} where worst-case convergence analysis is performed through the lens of control theory and finding Lyapunov functions. As such, the worst-case guarantees rely on smaller SDPs, but are only asymptotic.

\smallskip\noindent\textbf{Contributions.} Our main contributions are as follows
\begin{itemize}
    \item By extending the PEP framework to first-order block coordinate algorithms, we provide a unifying framework for the analysis of such algorithms. We illustrate the flexibility of our framework by analysing three variants of block coordinate algorithms: block coordinate descent, alternating minimization and a cyclic version of the random accelerated block coordinate descent algorithm from \cite{richtarik}.
    \item We compute the exact worst-case convergence rate of block coordinate algorithms over the class of smooth convex functions for a range of optimization setups. Furthermore, we provide sublinear upper and lower bound valid for the class of functions with coordinate-wise Lipschitz gradients, which is frequently considered when analysing such methods. For the block coordinate descent algorithm, our bounds outperforms the best known bound from \cite{Beck} by an order of magnitude. For the alternating minimization algorithm our bounds suggest that the bound from \cite{Beck} is asymptotically tight.
    \item We provide numerical worst-case convergence analysis for cyclic block coordinate algorithms under more natural and less restrictive assumptions than those usually made for the analysis of cyclic block coordinate algorithms \cite{Beck,Nesterov}.
    \item We show that acceleration schemes for random block-wise algorithms do not necessarily generalize to cyclic variants; we provide indeed a numerical lower bound on the convergence rate of the latter that suggests it is asymptotically worse than that of the original random version.
\end{itemize}
\textbf{Related work.} As mentioned above, cyclic block coordinate descent has attracted less attention than its random counterparts. Although a convergence rate have been established in \cite{Beck} for the unconstrained smooth convex minimization setting, it is quite conservative and the assumptions made on the first iterate are not very practical.
\newline
Relying on the performance estimation problem to establish worst-case convergence bounds for cyclic coordinate descent has been attempted in \cite{Shi}. Although the bound they obtain is significantly better then the bound in \cite{Beck}, it is established under very restrictive assumptions such as the equality of the dimensions of the blocks. In \cite{Bach}, an asymptotic worst-case convergence bound (not guaranteed to be tight) is established for random coordinate descent using Lyapunov functions.

\section{Block coordinate algorithms}
We consider the general unconstrained minimization setup
\begin{equation} 
    \min_{x \in \mathcal{R}^d} f(x),
\end{equation}
where the function $f$ is convex, differentiable and defined over the entire space $\mathcal{R}^d$. We assume that the optimal set $X^*$ is non-empty and select an arbitrary optimal point $x^{(*)} \in X^*$. As we will analyze block-coordinate algorithms, we further consider a partition of the space $\mathcal{R}^d$ into $p$ subspaces
\begin{equation}
    \mathcal{R}^d = \mathcal{R}^{d_1} \times ..... \times \mathcal{R}^{d_p},
\end{equation}
and introduce selection matrices $U_i\in \mathcal{R}^{n\times n_i}$ such that
\begin{equation*}
    \left( U_1,\dots,U_p \right) = \mathcal{I}_n,
\end{equation*}
allowing to write for every $x \in \mathcal{R}^d$
\begin{equation}\label{eqx_sep_p_blocks}
    x = (x_1;\dots;x_p) \; \text{with}  \;  x_i = U_i^{T}x \in \mathcal{R}^{d_i} \; \forall i \in 1,\dots,p.
\end{equation}
If $x$ is given as in \eqref{eqx_sep_p_blocks}, we also have that $ x= \sum_{i = 1}^p U_i x_i$.
\newline
\textbf{Definition 1.1} We define the partial gradient of $f$ in $x_i$ by
\begin{equation}
    \nabla_i f(x) \triangleq U_i^T \nabla f(x).
\end{equation}
We now present three specific cyclic block-coordinate algorithms that we will analyse, though our approach can be directly applied to a wide class of methods. First the cyclic block coordinate descent (CCD) performs at each iteration a gradient step with respect to a block of variables chosen in a cyclic order.
\begin{algorithm}[H]
\caption{Cyclic coordinate descent (CCD)}\label{algcap}
\begin{algorithmic}
\State \textbf{Input} starting point $x^{(0)} \in \mathcal{R}^d$ and step-size $\alpha$, number of cycles $K$, number of blocks $p$, $N = pK$.
\State For $n = 1 \dots N$ ,
\State \hspace{1.5cm} Set $i = \text{mod}(n,p) + 1$
\State \hspace{1.5cm} $x^{(n)} = x^{(n-1)} - \alpha U_i\nabla_if(x^{(n-1)})$
\State \textbf{Output} $x^{(N)}$ 
\end{algorithmic}
\end{algorithm}
Note that $K$ denotes the number of cycles, where each block of coordinate is updated once in a predefined order, and $p$ the number of blocks. Thus the number of partial gradient steps in Algorithm 1 is $N = pK$.
\newline
We also consider the cyclic alternating minimization (AM) where at each step we perform an exact minimization along a chosen block of coordinates instead of using a partial gradient step.
\begin{algorithm}[H]
\caption{Alternating minimization (AM)}\label{algcap}
\begin{algorithmic}
\State \textbf{Input} starting point $x^{(0)}\in \mathcal{R}^d$ and number of cycles $K$, number of blocks $p$, $N = pK$.
\State For $n = 1 \dots N$ ,
\State \hspace{1.5cm} Set $i = \text{mod}(n,p) + 1$
\State \hspace{1.5cm} $x^{(n)} = \arg \min_{z = x^{(n-1)}+U_i\Delta x_i, \; \Delta x_i \in \mathcal{R}^{d_i}} f(z)$
\State \textbf{Output} $x^{(N)}$ 
\end{algorithmic}
\end{algorithm}
Finally, we also analyse a deterministic cyclic version of the accelerated random coordinate descent algorithm from \cite{richtarik}, which we denote (CACD).

\begin{algorithm}[H]
\caption{Cyclic accelerated coordinate descent (CACD)}\label{algcap}
\begin{algorithmic}
\State \textbf{Input} starting point $z^{(0)} \in \mathcal{R}^d$, $\theta_{0} = \frac{1}{p}$, number of cycles $K$, number of blocks $p$, $N = pK$
\State For $n = 1 \dots N$ 
\State \hspace{1.5cm} Set $i = \text{mod}(n,p) + 1$
\State \hspace{1.5cm} $y^{(n-1)} = (1-\theta_k)x^{(n-1)} + \theta_{n-1} z^{(n-1)}$
\State \hspace{1.5cm} $z^{(n)} = z^{(n-1)} - \frac{1}{p \theta_n L} U_i\nabla_if(y^{(n-1)})$
\State \hspace{1.5cm} $x^{(n)} = y^{(n-1)} + p \theta_{n-1} (z^{(n)} -z^{(n-1)})$
\State \hspace{1.5cm} $\theta_{n} = \frac{\sqrt{\theta_{n-1}^4 + 4 \theta_{n-1}^2} -\theta_{n-1}^2}{2}$
\State \textbf{Output} $x^{(N)}$ 
\end{algorithmic}

\end{algorithm}

\section{The class of convex smooth functions and convex smooth interpolation}
In the sequel, We will consider two classes of functions in our worst-case analysis.
\newline
\newline
\textbf{Definition 1.2} $f$ is $L$-smooth if and only if
\begin{equation}
    \forall x \in \mathcal{R}^{d}, \; \forall h \in \mathcal{R}^{d},\; ||\nabla f(x+h)-\nabla f(x)|| \leqslant L ||h||
\end{equation}
where $||.||$ is the standard Euclidean norm. We denote the class of convex smooth functions by $\mathcal{F}_L$. However, in the context of block coordinate algorithms, most of the existing work assumes instead a form of coordinate-wise smoothness.
\newline
\newline
\textbf{Definition 1.3} Given a vector $\textbf{L} = (L_1,\dots,L_p)$ of $p$ nonnegative constants, $f$ is $\textbf{L}$-coordinate-wise smooth if and only if we have $\forall i \in 1,\dots,p$
\begin{equation}
    \forall x \in \mathcal{R}^{d}, \; \forall h_i \in \mathcal{R}^{d_i},\; ||\nabla_i f(x+U_ih_i)-\nabla_i f(x)|| \leqslant L_i ||h_i|| 
\end{equation}
We denote the class of convex coordinate-wise smooth functions by $\mathcal{F}^{\text{coord}}_{\textbf{L}}$.
\newline
\newline
We will focus on the analysis of $L$-smooth function $\mathcal{F}_L$, which proves simpler. However, our results will have direct implications for the class of coordinate-wise smooth functions $\mathcal{F}^{\text{coord}}_L$ thanks to the following lemma, whose proof is direct.
\newline
\textbf{Lemma 1.1.} Given any vector $\textbf{L} \in \mathcal{R}^{p}_{+}$, we have $\mathcal{F}_{L_{min}} \subset \mathcal{F}^{\text{coord}}_{L} \subset \mathcal{F}_{\overline{L}}$ where $\overline{L} = \sum^p_{i = 1} L_i$ and $L_{min} = \min_{i = 1,\dots,p}L_i$.
\newline
\newline
This lemma shows that our exact worst-case bound for block coordinate algorithms on (globally) smooth functions also automatically provide valid lower and upper bounds for the class of coordinate-wise smooth functions.
\section{Smooth convex interpolation.}
We now review the smooth convex interpolation result of \cite{Taylor}, which is crucial when deriving exact worst-case bounds for first-order optimization algorithms over the class of smooth convex functions, and will be instrumental in our analysis in Section \ref{PEP}.
\newline
\newline
\textbf{Definition 2.2.} A set $\{(x^{(n)},g^{(n)},f^{(n)})\}_{n = 1,\dots N} \subset \mathcal{R}^{d} \times \mathcal{R}^{d} \times \mathcal{R}$ such as $x^{(n)}_i = U_i^T x^{(n)}$ and $g^{(n)}_i = U_i^T g^{(n)}$, $\forall i = 1,\dots,p$ is $\mathcal{F}_L$-interpolable if and only if there exists a function $f \in \mathcal{F}_L$ such that
\begin{equation}
    \begin{aligned}
          & f^{(n)} = f(x^{(n)}) \; \forall n \in \{1,\dots,N\}\\
          &  g^{(n)}_{i} = \nabla_i f(x^{(n)}) \; \forall i \in \{1,\dots,p\}, \; \forall n \in \{1,\dots,N\}
    \end{aligned}
\end{equation}
\newline
\textbf{Theorem 2.2.} The set $\{(x^{(n)},g^{(n)},f^{(n)})\}_{n = 1,\dots,N}$ is $\mathcal{F}_L$-interpolable if and only if
\begin{multline}\label{interp conds}
    \forall n,l \in \{1,\dots,N\},\\ f^{(n)} \geqslant f^{(l)} + \sum^p_{i = 1} \langle g^{(l)}_{i},x^{(n)}_{i}-x^{(l)}_{i}\rangle + \frac{1}{2L} \sum^p_{i=1}||g^{(n)}_{i}-g^{(l)}_{i}||^2
\end{multline}
Any set $\{(x^{(n)},g^{(n)},f^{(n)})\}_{n = 1,\dots N}$ satisfying $\eqref{interp conds}$ is thus consistent with an actual globally defined smooth convex function. We refer the reader to \cite{Taylor} for a detailed proof of this result expressed in terms of the full gradient of $f$. The reformulation in terms of partial gradients shown in \eqref{interp conds} is direct and more suitable to handle block coordinate algorithms.
\section{Worst-case behaviour of \newline coordinate descent-like algorithms.} \label{PEP}
We now present the performance estimation framework that will allow us to derive exact worst-case bounds for block coordinate algorithms. The idea behind the PEP framework is to cast the performance analysis of an optimization method as an optimization problem itself over the class of functions $\mathcal{F}_L$ (see \cite{Drori} and \cite{Taylor}).
\newline
\newline
For pedagogical reasons, we start by giving a simple example of a PEP for 2-block coordinate descent performing only one cycle. This can be written $x^{(1)} = x^{(0)} - \alpha U_1\nabla_1f(x^{(0)})$ and $x^{(2)} = x^{(1)} - \alpha U_2\nabla_2f(x^{(1)})$ . We choose as a performance measure the difference between the objective function value after one cycle and the optimal value. We also assume that there exists a minimizer $x^{(*)}$ of $f$ such that for a given fixed positive constant $R$, we have that $||x^{(0)}-x^{(*)}||^2 \leqslant R^2$ i.e $R$ is a bound on the distance from the starting point $x^{(0)}$ to a minimizer.
In this simple case, we are looking for the $L$-smooth convex function that performs worst after one cycle, meaning the function that maximizes our performance measure $f(x^{(2)})- f(x^{(*)})$. The PEP can then be written as follow
\begin{equation}\label{eqconceptual PEP}
    \begin{aligned}
      \mathcal{W}_L(p = 2,K = 1,f,R)  = & \max f(x^{(2)}) - f(x^{(*)})\\
      & f \in \mathcal{F}_{L} \\
      & ||x^{(0)}-x^{(*)}||^2 \leqslant R^2 \\
      &  x^{(1)} = x^{(0)} - \alpha U_1\nabla_1f(x^{(0)})\\
      &  x^{(2)} = x^{(1)} - \alpha U_2\nabla_2f(x^{(1)})\\
      & x^{(*)} \; \text{is a minimizer of $f$}
    \end{aligned}
\end{equation}
\newline
\newline
Problem $\eqref{eqconceptual PEP}$ is an infinite-dimensional problem over the class of functions $\mathcal{F}_L$ and cannot be directly solved numerically. Convex interpolation transforms such problems into tractable finite-dimensional convex semidefinite programs (SDPs). The idea is to replace $f$ by a set of variables of the form $\{(x^{(n)},g^{(n)},f^{(n)}))\}_{n = 0,1,2,*}$ that satisfy the interpolation conditions $\eqref{interp conds}$
\begin{equation}\label{eqPEPrelaxation}
    \begin{aligned}
      \mathcal{W}_L(p =2,K = 1,R)  = & \max f^{(2)} - f^*\\
      & \{(x^{(n)},g^{(n)},f^{(n)}))\}_{n = 0,1,2,*} \; \text{is} \\& \mathcal{F}_L- \text{interpolable}\\
      & ||x^{(0)} - x^{(*)}||^2 \leqslant R^2 \\
      &  x^{(1)}_1 = x^{(0)}_1 - \alpha g^{(0)}_1\\
      & x^{(1)}_2 = x^{(0)}_2\\
      & x^{(2)}_1 = x^{(1)}_1\\
      &  x^{(2)}_2 = x^{(1)}_2 - \alpha g^{(1)}_2\\
      & g^{(*)}_1 = g^{(*)}_2 = 0
    \end{aligned}
\end{equation}
Problem \eqref{eqPEPrelaxation} is indeed an exact reformulation of problem \eqref{eqconceptual PEP} because the interpolation conditions \eqref{interp conds} guarantee the existence of a smooth convex function that interpolates the set of variables in problem \eqref{eqPEPrelaxation} and therefore provides a solution of problem \eqref{eqconceptual PEP}.
We now give a more general form of PEPs for block coordinate algorithms with a fixed arbitrary number of blocks and cycles
\begin{equation}\label{eqPEP relaxation}
    \begin{aligned}
      \mathcal{W}_L(p,K,R)  = & \max \mathcal{P}(\mathcal{O}_n)\\
      & \mathcal{O}_n \; \text{satisfies \eqref{interp conds}}\\
      & \mathcal{I}(\mathcal{O}_n) \leqslant R^2 \\
      & \{x^{(n)}\}_{n =1,\dots,N} \; \text{are computed}\\
      & \text{by the considered algorithm}
    \end{aligned}
\end{equation}
where $\mathcal{O}_n$ denotes the set $\{(x^{(n)},g^{(n)},f^{(n)}))\}_{n = 0,\dots,N,*}$, $\mathcal{P}(\mathcal{O}_n)$ denotes a performance criterion and $\mathcal{I}(\mathcal{O}_n)$ a measure of the optimality of the initial iterate such as $||x^{(0)}-x^{(*)}||^2$.
\newline
\newline
For the (CCD) and (CACD) algorithm, all constraints in the PEP are linear equalities involving the iterates $x^{(n)}$ and the gradients $g^{(n)}$. This can be written in term of Gram matrices of vectors $\{x^{(n)}\}$ and $\{g^{(n)}\}$. Indeed, the scalar products in $\eqref{interp conds}$ can then be expressed thanks to $x^{(n)T} g^{(l)} = \sum_{i=1}^p x^{(n)}_ig^{(l)}_i$. For the (AM) algorithm, the partial minimization along a coordinate amounts to imposing a partial gradient equal to zero along this coordinate, which can also be written in terms of Gram matrices. We refer the reader to \cite{Taylor} for the detailed procedure to transform a PEP into a SDP in the case of full gradient methods. The procedure we use to transform $\eqref{eqPEP relaxation}$ into a SDP is very similar to one presented in \cite{Taylor}. Note that this procedure also requires that $\mathcal{P}$ and $\mathcal{I}$ can be written in terms of gram matrices, which is the case for the choices of $\mathcal{P}$ and $\mathcal{I}$ that we make in the sequel. The main difference is that we have to consider multiple blocks by defining a distinct Gram matrix for each of them. For a given number $p$ of blocks, we solve an SDP involving $p$ semidefinite matrices.
\newline
\newline
\textbf{Theorem 3.1.} If $f$ is a convex $L$-smooth function, then the performance of $N$ steps of  Algorithms 1, 2 or 3 over $f$, denoted by $\mathcal{P}(f)$, verifies the following
\begin{equation}
\label{eq PEP bound}
    \mathcal{P}(f) \leqslant \mathcal{W}_L(p,K,R)
\end{equation}
\newline
\textbf{Proof.} This is a direct consequence of the optimization problem \eqref{eqPEP relaxation} that defines $\mathcal{W}_L(p,K,R)$.
\newline
\newline
The previous theorems allow us to analyze the exact performance of block coordinate algorithms on $L$-smooth functions. However, we already noted that rates are in general expressed over the class of \emph{coordinate-wise smooth functions}. The next result allows us to derive bounds for this class of functions.
\newline
\newline
\textbf{Theorem 3.2} For any given vector of $p$ nonnegative constants $\textbf{L} = (L_1,\dots,L_p)$ and any block coordinate algorithm, we denote by $\mathcal{W}^{\text{coord}}_{\textbf{L}}(p,K,R)$ the worst-case of some given algorithm over the class of coordinate-wise smooth functions $\mathcal{F}^{\text{coord}}_{\textbf{L}}$ and we have
\begin{multline}
 \mathcal{W}_{L_{min}}(p,K,R) \leqslant \mathcal{W}^{\text{coord}}_{\textbf{L}}(p,K,R) \leqslant \mathcal{W}_{\overline{L}}(p,K,R)
\end{multline}
where $\mathcal{W}_{L_{min}}(p,K,R)$, $\mathcal{W}_{\overline{L}}(p,K,R)$ denote respectively the bound given by the PEP for the classes of functions $\mathcal{F}_{L_{min}}$ and $\mathcal{F}_{\overline{L}}$.
\newline
\newline
\textbf{Proof.} Let us consider $f \in \mathcal{F}^{\text{coord}}_{\textbf{L}}$. Thanks to Lemma 1.1, we know that $f \in \mathcal{F}_{\overline{L}}$ thus the performance of the algorithm $\mathcal{P}(f)$ on $f$ verifies
\begin{equation*}
    \mathcal{P}(f) \leqslant \mathcal{W}_{\overline{L}}(p,K,R)
\end{equation*}
Since this is true for every function $f \in \mathcal{F}^{\text{coord}}_{\textbf{L}}$, we have that
\begin{equation*}
    \mathcal{W}^{\text{coord}}_{\textbf{L}}(p,K,R) \leqslant \mathcal{W}_{\overline{L}}(p,K,R)
\end{equation*}
\newline
Similarly by considering the other inclusion given in Lemma 1.1, we obtain the other inequality of Theorem 3.2.
\newline
\newline
\textbf{Remark.}Theorem 3.2 applies when the same algorithm, including its coefficient value, is considered for the three classes. 
Its use requires some care, as algorithm parameters are frequently made implicitly dependent on the function class parameters; step-sizes are e.g. typically described as $\frac{h}{L}$. The application of Theorem 3.2 in such cases require thus first fixing the values of these parameters independently of $L$.
\newline
We will focus on particular instances of PEPs by choosing different performance criteria $\mathcal{P}$ and starting iterate conditions $\mathcal{I}$. For the performance criterion, we will consider the difference between the function value at the last iterate and the optimal value of the function
\begin{equation*}
    \mathcal{P}(\mathcal{O}_n) = f^{(N)} - f^*
\end{equation*}
and the squared gradient norm of the last iterate
\begin{equation*}
    \mathcal{P}(\mathcal{O}_n) = ||g^{(N)}||^2
\end{equation*}
As starting iterate conditions we consider the following two assumptions. First, the usual setting used for cyclic block coordinate descent presented in \cite{Beck}:
\newline
\newline
\textbf{Setting ALL.} We consider the following assumption on all the iterates
\begin{equation}
    \max_{x^{(*)}\in x^{(*)}} \max_{k=1\dots,K} ||x^{(pk)}-x^{(*)}|| \leqslant R_a
\end{equation}
\newline
Setting ALL is an adaptation for the PEP framework of the assumptions usually made for the analysis of cyclic coordinate algorithms. In \cite{Beck,Nesterov}, it is indeed assumed that the set $S = \{x\in \mathcal{R}^d:\; f(x)\leqslant f(x^{(0)})\}$ is compact which implies that
\begin{equation}
    ||x^{(pk)} -x^{(*)}|| \leqslant R(x^{(0)}), \; \forall k \in 1,\dots,K
\end{equation}
with $R(x^{(0)})$ defined as
\begin{equation}
    R(x^{(0)}) = \max_{x \in x^{(*)}} \max_{x \in \mathcal{R}^{d}} \{||x-x^{(*)}||:\; f(x) \leqslant f(x^{(0)})\}
\end{equation}
\textbf{Theorem 3.2} Let $\{x^{(n)}\}_{n = 1\dots,N}$ be a sequence generated by the (CCD) algorithm with a constant step-size $\alpha \leqslant \frac{1}{L}$ then
\begin{equation}
    f(x^{(N)})-f^* \leqslant \frac{4}{\alpha} (1+p\alpha^2L^2)\frac{p}{N+8}R_a^2
\end{equation}
Theorem 3.2 is adapted from [Theorem 3.6] in \cite{Beck} whose proof remains valid when $R(x^{(0)})$ is replaced by $R_a$.
\newline
\newline
Though theoretically convenient, Setting ALL. is not very natural and may be difficult to verify in practice. In addition, Setting ALL. can prove to be unusable for certain class of functions. Indeed, consider the family of smooth functions $f_{\epsilon}(x,y) =(x-y)^2 + \epsilon (x^2+y^2)$ and the initial point $(x_0 =1,y_0=-1)$ for 2-block coordinate algorithm. We have that $R(x^{(0)}) = R_a = \frac{1}{\sqrt{\epsilon}}$ which tends to infinity when $\epsilon$ tends to zero. For $\epsilon$ small enough the bounds obtained in this setting are very conservative and do not give useful information about the performance of the algorithm. Therefore, we will also consider a more classical setting, albeit less frequently used in the context of deterministic block coordinate algorithms:
\newline
\newline
\textbf{Setting INIT.} Given the starting point of the block coordinate algorithm $x^{(0)}$ and an optimal point of the function $x^{(*)}$, we have that
\begin{equation}\label{eqnew R}
     ||x^{(0)}-x^{(*)}||^2 \leqslant R_i^2,
\end{equation}

\section{Bounds on the worst-case of coordinate descent algorithms using the PEP framework}\label{exps}
\begin{figure}[H]
    \centering
        \includegraphics[width=0.4\textwidth]{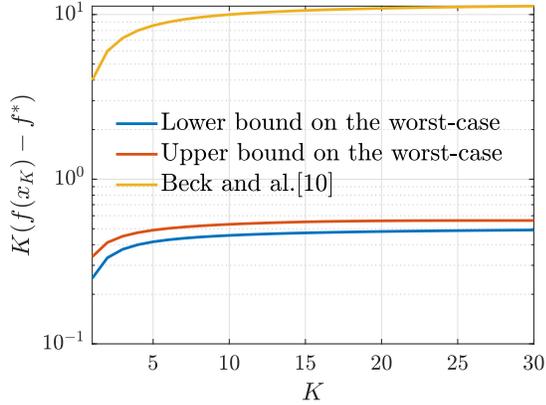}
    \caption{Comparison between the PEP bounds and the theoretical bound in \cite[Theorem 3.6]{Beck} multiplied by the number of cycles $K$, for the 2-block (CCD) in the Setting ALL.}
    \label{fig 1.}
\end{figure}

\begin{figure}[H]
    \centering
        \includegraphics[width=0.4\textwidth]{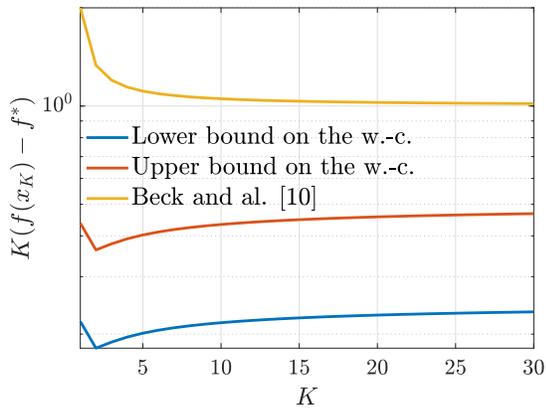}
    \caption{Comparison between the PEP bounds and the theoretical bound in \cite[Theorem 5.2]{Beck} multiplied by the number of cycles $K$, for the 2-block (AM) algorithm in Setting ALL.}
    \label{fig 2.}
\end{figure}

\begin{figure}[H]
    \centering
        \includegraphics[width=0.4\textwidth]{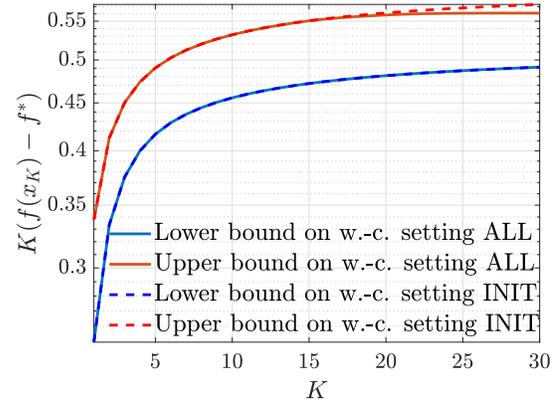}
    \caption{Comparison between the PEP upper bounds multiplied by the number of cycles $K$, for the 2-block (CCD) in Settings ALL (full lines) and INIT.(dashed lines)}
    \label{fig 3.}
\end{figure}

\begin{figure}[H]
    \centering
        \includegraphics[width=0.4\textwidth]{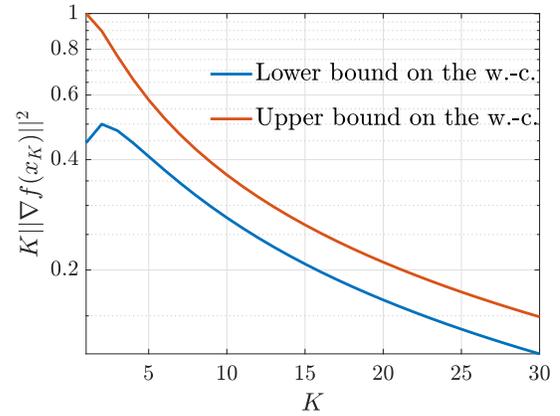}
    \caption{PEP bounds multiplied by the number of cycles $K$ on the residual gradient norm for the 2-block (CCD) in Setting all.}
    \label{fig 4.}
\end{figure}

\begin{figure}[H]
    \centering
        \includegraphics[width=0.4\textwidth]{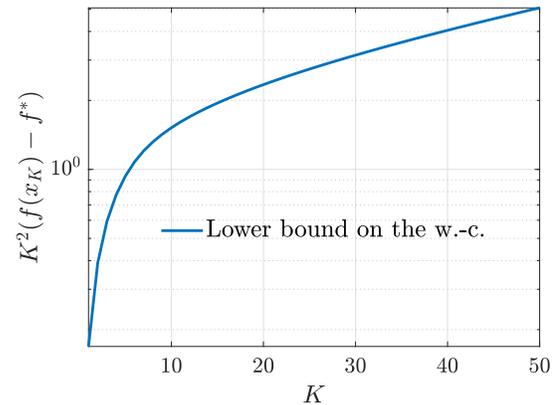}
    \caption{PEP lower bound multiplied by the number of cycles squared $K^2$, for the 2-block (CACD), which indicates that the convergence rate is slower than $\mathcal{O}(\frac{1}{K^2})$.}
    \label{fig 5.}
\end{figure}

We now exploit our PEP-based approach to revisit and in some cases significantly improve the bounds on the block coordinate algorithms defined in the first section. For simplicity only, we focus on algorithms on the case of two blocks $(p = 2)$. Furthermore, since $\mathcal{W}$ is proportional to $L$ for step-sizes $\alpha = \frac{h}{L}$, we consider without loss of generality smoothness constants $L_1 = L_2 = 1$. For the same reason we choose $R_a = R_i =1$. Since a $\mathcal{O}(\frac{1}{K})$ rate of convergence is expected in most cases, we show the evolution of our performance criterion multiplied by $K$ to facilitate the analysis.
\newline
\newline
In Figure $\ref{fig 1.}$, we provide the upper and lower bound defined in Theorem 3.2 for the worst-case of the 2-block (CCD) with a constant step-size $\alpha = \frac{1}{2}$. Our bounds are sublinear and improve by one order of magnitude the best know theoretical bound derived in \cite[Theorem 3.6]{Beck}. The experiments presented in Figure $\ref{fig 1.}$ are performed under Setting ALL defined in the previous section. In Figure \ref{fig 2.} we also provide upper and lower bounds for the alternating minimization algorithm in setting ALL and compare it to the bound provided in \cite[Theorem 5.2]{Beck}. Figure 2 shows that our PEP-based approach improves bounds by a factor of two.
\newline
\newline
In order to illustrate the flexibility and usefulness of our framework, we provide bounds for the 2-block (CCD) algorithm in the Setting INIT with two different performance criteria: the objective accuracy in Figure $\ref{fig 3.}$ and the squared gradient norm of the last iterate in Figure $\ref{fig 4.}$. In Figure $\ref{fig 3.}$ we compare the upper bounds obtained for both settings for block coordinate descent. Our results show that the convergence of cyclic block coordinate descent can also be established under the weaker assumption of Setting INIT, and the performance results suggest that the stronger assumptions made in setting ALL do not yield a significant improvement in terms of performance. Figure \ref{fig 4.} suggests that the squared residual gradient norm converges faster then $\mathcal{O}(\frac{1}{K})$.
\newline
\newline
Finally, we provide an exact worst-case bound for the cyclic 2-block version of the random accelerated coordinate algorithm (CACD) derived in \cite{richtarik} over the class of $1$-smooth functions which, thanks to Theorem 3.2, also gives a lower bound on the worst-case for the class of function $\mathcal{F}^{\text{coord}}_L$. The random version of (CACD) from \cite{richtarik} has a $\mathcal{O}(\frac{1}{K^2})$ rate of convergence. We do not observe acceleration in Figure $\ref{fig 5.}$ in the sense that the rate of convergence is slower than $ \mathcal{O} (\frac{1}{K^2})$. This indicates that using randomness in the choice of the block of coordinates to update plays a crucial role in the acceleration of block coordinate algorithms. To investigate this further, we adapted our PEP framework to compare the cyclic and random versions of the algorithm presented in \cite{richtarik}. The PEP framework is usually only able to handle deterministic algorithms. To circumvent this issue, we write a PEP that computes simultaneously all possible choices of coordinate steps, and use as a performance criterion the worst-case average of the performances of each combination, which corresponds exactly to the expectation of the performance of the random accelerated coordinate descent in \cite{richtarik}. For readability, we refer to this approach as the worst-case of the average, and we compare it to the worst-case of the deterministic version of the algorithm in \cite{richtarik} for each possible combination of steps. We  consider again 2-blocks of coordinates for $N$ partial steps of gradient. Since there are $2^N$ possible combinations of steps and the dimension of the SDP grows with that number, we only present preliminary results for $N = 4$. The worst-case of the average in this case denoted by $\mathcal{W}_4$ is equal to $\mathcal{W}_4 = 0.1046$.
\newline
\newline
Table $\ref{tabtable1}$ give the worst-case for each possible combination of $N = 4$ partial steps. Note that all the worst-cases in Table $\ref{tabtable1}$  are larger than $\mathcal{W}_4$ which tends to confirm that the usual acceleration schemes used for coordinate descent are specific to random coordinate descent. These preliminary results tend to indicate that for general smooth functions the best deterministic choice of coordinates is the cyclic one.
\begin{table}[H]
    \centering
    \begin{tabular}{|c||c|}
        \hline
        Ordered choice of steps  & Worst-case \\ 
        \hline
        1    1    1   1 \newline  2    2    2    2    & 0.5\\ 
        \hline
        1    1    1    2 \newline 2    2    2    1  & 0.25517 \\
        \hline
        2    2    1    1 \newline 1    1    2    2 & 0.23462 \\
        \hline
        2    1    1    1 \newline  1    2    2    2 & 0.19905 \\
        \hline
        1    1    2    1 \newline  2    2    1    2 & 0.19574 \\
        \hline
               1    2    1    1 \newline 2    1    2    2 &  0.16453 \\
        \hline
         1    2    2    1  \newline  2    1    1    2   & 0.14988 \\
        \hline
        2    1    2    1 \newline  1    2    1    2&  0.14429 \\
        \hline
         $\mathcal{W}_4$ &    0.1046 \\   
        \hline
\end{tabular}
    \caption{Worst-cases of the cyclic 2-blocks version of the accelerated coordinate descent algorithm presented in \cite{richtarik} for $N = 4$. The first column indicates the index of the coordinate chosen for each step and $\mathcal{W}_4$ denotes the worst-case of the average defined in section \ref{exps} }
    \label{tabtable1}
\end{table}


\section{Conclusion.}
We have developed a flexible framework for the automated worst-case analysis of block coordinate algorithms, and provided exact worst-case bounds over the class of smooth functions, which lead to upper and lower bounds on the class of coordinate-wise smooth functions. We have provided improved numerical bounds, sometimes by an order of magnitude, for three types of block coordinate algorithms:  block coordinate descent (CCD), alternating minimization (AM) and a cyclic version of the accelerated random coordinate descent in \cite{richtarik} (CACD). In addition, we highlighted the importance of randomness for existing acceleration schemes, since our numerical experiments suggest that deterministic cyclic algorithms do not accelerate i.e they do not achieve a $\mathcal{O}(\frac{1}{K^2})$ rate of convergence. Further research could involve developing interpolation conditions for the class of coordinate-smooth functions, performing more numerical experiments in a wider range of settings, with more blocks and different step-sizes, as well as searching with our PEP-based approach an efficient acceleration schemes for cyclic block coordinate algorithms over the class of coordinate-wise smooth functions.

\section*{Acknowledgement}
Y. Kamri is supported by the European Union’s MARIE SKŁODOWSKA-CURIE Actions Innovative Training Network (ITN)-ID 861137, TraDE-OPT.

\bibliography{biblio.bib}

\end{document}